\documentclass[12pt]{amsart}
\usepackage{tikz}
\usetikzlibrary{decorations.pathmorphing}
\title{On the enumeration of Hoppy's walks}

\author{Helmut Prodinger}
\address{Helmut Prodinger\\
	Department of Mathematical Sciences\\
	Stellenbosch University\\
	7602 Stellenbosch\\
	South Africa}
\email{hproding@sun.ac.za}

\begin{document}

	\maketitle
	
	\section{Hoppy walks}
	
	Deng and Mansour \cite{deng} introduce a rabbit named Hoppy and let him move according to certain rules.
	At that stage, we don't need to know the rules. Eventually, the enumeration problem is one about $k$-Dyck paths.
	The up-steps are $(1,k)$ and the down-steps are $(1,-1)$.
	
	
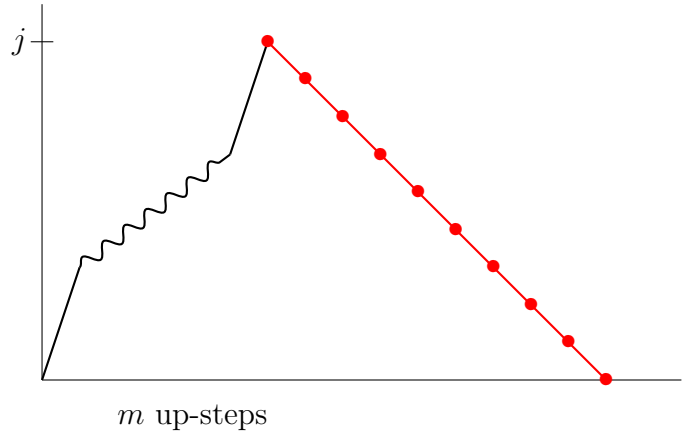
\begin{center}
	\begin{tikzpicture}[scale=0.5]
	\draw (0,0) -- (17,0);
	\draw (0,0) -- (0,10);
	\draw[thick](0,0)--(1,3);
	\draw[thick](5,6)--(6,9);
	\draw [decorate,decoration=snake,thick] (1,3) -- (5,6);
	\foreach \i in {0,...,8}
	{\draw[thick,red](6+\i,9-\i)--(7+\i,8-\i);
		\node[thick, red] at (6+\i,9-\i){$\bullet$};
		}
	\node[thick, red] at (6+5+4,9-5-4){$\bullet$};

	\node at (4,-1){$m$ up-steps};
	 
	\node[thick] at (-1+0.4,9){$j$};

	\draw (-0.3,9) --(0.3,9);
	\end{tikzpicture}
	\end{center}	
	The question is about the length of the sequence of down-steps printed in red. Or, phrased differently, how many $k$-Dyck paths end on level $j$, after $m$ up-steps, the last step being an up-step. The recent paper \cite{jcmcc} contains similar computations, although without the restriction that the last step must be an up-step. 
	
	Counting the number of up-steps is enough, since in total, there are $m+km=(k+1)m$ steps. The original description of Deng and Mansour is a reflection of this picture, with up-steps of size 1 and down-steps of sice $-k$, but we prefer it as given here, since we are going to use the adding-a-new-slice method, see \cite{FP, Prodinger-handbook}. A slice is here a run of down-steps, followed by an up-step. 
	The first up-step is treated separately, and then $m-1$ new slices are added. We keep track of the level after each slice, using a variable $u$. The variable $z$ is used to count the number of up-steps.
	
	Deng and Mansour work out a formula which comprises $O(m)$ terms. Our method leads only to a sum of $O(j)$ terms.
	
	The following substitution is essential for adding a new slice:
	\begin{equation*}
u^j\longrightarrow z\sum_{0\le h \le j} u^{h+k}=\frac{zu^k}{1-u}(1-u^{j+1}).
	\end{equation*}
	Now let $F_m(z,u)$ we the generating function according to $m$ runs of down-steps. The substitution leads to
	\begin{equation*}
F_{m+1}(z,u)=\frac{zu^k}{1-u}F_m(z,1)-\frac{zu^{k+1}}{1-u}F_m(z,u),\quad F_0(z,u)=zu^k.
	\end{equation*}
	Let $F=\sum_{m\ge0}F_m$, then
	\begin{equation*}
F(z,u)=zu^k+\frac{zu^k}{1-u}F(z,1)-\frac{zu^{k+1}}{1-u}F(z,u),
	\end{equation*}
	or
	\begin{equation*}
		F(z,u)\frac{1-u+zu^{k+1}}{1-u}=zu^k+\frac{zu^k}{1-u}F(z,1).
	\end{equation*}	
	The equation $1-u+zu^{k+1}=0$ is famous when enumerating $(k+1)$-ary trees. Its relevant combinatorial solution (also the only one being analytic at the origin) is
		\begin{equation*}
		\overline{u}=\sum_{\ell\ge0}\frac1{1+\ell(k+1)}\binom{1+\ell(k+1)}{\ell}z^\ell.
	\end{equation*}
	Since $u-\overline{u}$ is a factor of the LHS, is must also be a factor of the RHS, and we can compute (by dividing out the factor
	$(u-\overline{u})$) that
	\begin{equation*}
\frac{zu^k(1-u+F(z,1))}{u-\overline{u}}=-zu^k.
	\end{equation*}
Thus
\begin{equation*}
F(z,u)=zu^k\frac{\overline{u}-u}{1-u+zu^{k+1}}.
\end{equation*}	


The first factor has even a combinatorial interpretation, as a description of the first step of the path. It is also clear from this that the level reached is $\ge k$ after each slice. We don't care about the factor $zu^k$ anymore, as it produces only a simple shift. The main interest is now how to get to the coefficients of
\begin{equation*}
 \frac{\overline{u}-u}{1-u+zu^{k+1}}
\end{equation*}	
in an efficient way. There is also the formula
\begin{equation*}
1-u+zu^{k+1}=(\overline{u}-u)\Big(1-z\frac{u^{k+1}-\overline{u}^{k+1}}{u-\overline{u}}\Big),
\end{equation*}
but it does not seem to be useful here.

First we deal with the denominators
	\begin{equation*}
S_j:=[u^j]\frac1{1-u+zu^{k+1}}=\sum_{0\le m\le j/k}(-1)^m\binom{j-km}{m}z^m.
	\end{equation*}
	One way to see this formula is to prove by induction that the sums $S_j$ satisfy the recursion
	\begin{equation*}
S_j-S_{j-1}+zS_{j-k-1}=0
	\end{equation*}
	and initial conditions $S_0=\dots =S_k=1$. In \cite{jcmcc} such expressions also appear as determinants.
	Summarizing,	
	\begin{equation*}
\frac1{1-u+zu^{k+1}}=\sum_{m\ge0}(-1)^mz^m\sum_{j\ge km }\binom{j-km}{m}u^j.
	\end{equation*}
	Now we read off coefficients:
\begin{multline*}
		[u^j]\frac{\overline{u}}{1-u+zu^{k+1}}\\=\sum_{0\le m\le j/k}(-1)^m\binom{j-km}{m}z^m
		\sum_{\ell\ge0}\frac1{1+\ell(k+1)}\binom{1+\ell(k+1)}{\ell}z^\ell
		\end{multline*}
and further
\begin{multline*}
[z^n][u^j]\frac{\overline{u}}{1-u+zu^{k+1}}\\=\sum_{0\le m\le j/k}(-1)^m\binom{j-km}{m}
\frac1{1+(n-m)(k+1)}\binom{1+(n-m)(k+1)}{n-m}.
\end{multline*}
The final answer to the Deng-Mansour enumeration (without the shift) is
\begin{multline*}
\sum_{0\le m\le j/k}(-1)^m\binom{j-km}{m}
\frac1{1+(n-m)(k+1)}\binom{1+(n-m)(k+1)}{n-m}\\
-(-1)^n\binom{j-1-kn}{n}.
\end{multline*}
If one wants to take care of the factor $zu^k$ as well, one needs to do the replacements $n\to n+1$ and $j\to j+k$ in the formula just derived. That enumerates then the $k$-Dyck paths ending at level $j$ after $n$ up-steps, where the last step is an up-step.

\section{An application}

The encyclopedia of integer sequences \cite{OEIS} has the sequences A334680, A334682, A334719, (with a reference to \cite{AHS}) which is the total number of down-steps of the last down-run, for $k=2,3,4$. So, if the path ends on level $j$, the contribution to the total is $j$. 

All we have to do here is to differentiate
\begin{equation*}
F(z,u)=zu^k\frac{\overline{u}-u}{1-u+zu^{k+1}}.
\end{equation*}	
w.r.t. $u$, and then replace $u$ by 1. The result is
\begin{equation*}
\frac{\overline{u}}z-\overline{u}-\frac1z,
\end{equation*}
and the coefficient of $z^m$ therein is
\begin{align*}
\frac1{1+(m+1)(k+1)}\binom{1+(m+1)(k+1)}{m+1}-\frac1{1+m(k+1)}\binom{1+m(k+1)}{m}.
\end{align*}
I don't know how this was derived in \cite{AHS}, but it is more fun to figure out things for oneself!

We hope to report about more applications soon.

\section{Hoppy's early adventures}

Now we investigate what Hoppy does after his first up-step; he might follow with $0,1,\dots,k$ down-steps. Eventually, we want to sum all these steps (red in the picture).
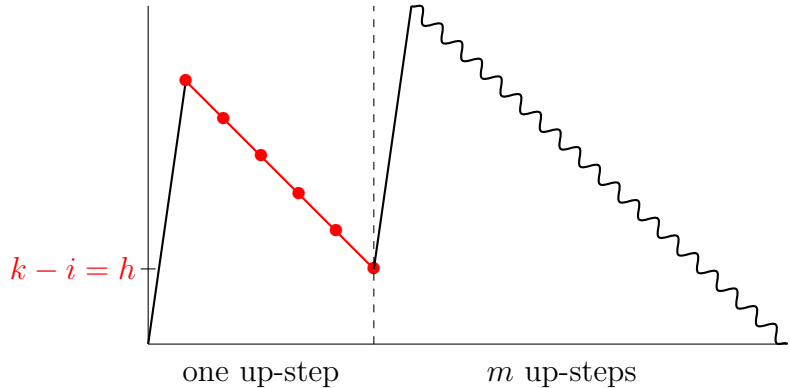
\begin{center}
	\begin{tikzpicture}[scale=0.5]
	\draw (0,0) -- (17,0);
	\draw (0,0) -- (0,9);
	\draw[thick](0,0)--(1,7);
	\foreach \i in {0,...,4}
	{\draw[thick,red](1+\i,7-\i)--(2+\i,6-\i);
		\node[thick, red] at (1+\i,7-\i){$\bullet$};
	}
	\node[thick, red] at (6,2){$\bullet$};
	\draw[thick](6,2)--(7,9);
	\draw [decorate,decoration=snake,thick] (7,9) -- (17,0);
	
	\draw(-0.2,2)--(0.2,2);
	\node[thick,red] at (-2,2){$k-i=h$};
\draw[ dashed](6,0)--(6,9);
\node[ ] at (3,-0.8){one up-step};
\node[ ] at (11,-0.8){$m$ up-steps};
	\end{tikzpicture}
\end{center}

A new slice is now an up-step, followed by a sequence of down-steps. The substitution of interest is:
 \begin{equation*}
	u^i\rightarrow  z\sum_{0\le h\le i+k} u^h=\frac{z}{1-u}-\frac{zu^{i+k+1}}{1-u}.
\end{equation*}
Furthermore
\begin{equation*}
	F_{h+1}(z,u)=\frac{z}{1-u}F_h(z,1)-\frac{zu^{k+1}}{1-u}F_h(z,u),
\end{equation*}
and $F_0=u^h$, the starting level.

We have
\begin{equation*}
H(z,u)=\sum_{h\ge0}F_h(z,u)=u^h+\frac{z}{1-u}H(z,1)-\frac{zu^{k+1}}{1-u}H(z,u)
\end{equation*}
or
\begin{equation*}
	H(z,u)(1-u+zu^{k+1}) =u^h(1-u)+zH(z,1)
\end{equation*}
Plugging in $\overline{u}$ into the RHS gives 0:
\begin{equation*}
zH(z,1)=-\overline{u}^h(1-\overline{u}),
\end{equation*}
and
\begin{equation*}
	H(z,u) =\frac{u^h(1-u)-\overline{u}^h(1-\overline{u})}{1-u+zu^{k+1}}.
\end{equation*}
But we only need $H(z,0)$, since we return to the $x$-axis at the end:
  \begin{equation*}
  	H(z,0) = [h=0]+\overline{u}^{h+1}-\overline{u}^h  .
  \end{equation*}
 The total contribution of red steps is then
 \begin{equation*}
k+\sum _{h=0}^k(k-h)(\overline{u}^{h+1}-\overline{u}^h)=\sum _{h=1}^k\overline{u}^{h};
 \end{equation*}
the coefficient of $z^m$ in this is the total contribution. 
Since $\overline{u}=1+z\overline{u}^{k+1}$,
there is the further simplification
\begin{equation*}
-1+\frac1z+\frac1{1-\overline{u}}=\sum_{m\ge1}\frac{k}{m+1}\binom{(k+1)m}{m}z^m.
\end{equation*}
The proof of this is as follows. Let $m\ge1$, then
\begin{align*}
[z^m]	\Bigl(-1+\frac1z+\frac1{1-\overline{u}}\Bigr)&=-[z^m]\frac1{z\overline{u}^{k+1}}\\
&=-[z^{m+1}]\sum_{\ell\ge0}\frac {-(k+1)}{(k+1)\ell -(k+1)}\binom{(k+1)\ell -(k+1)}{\ell}z^\ell\\
&=[z^{m+1}]\sum_{\ell\ge0}\frac {(k+1)}{(k+1)(\ell-1)}\binom{(k+1)(\ell-1)}{\ell}z^\ell\\
&= \frac {(k+1)}{(k+1)m}\binom{(k+1)m}{m+1}=\frac{k}{m+1}\binom{(k+1)m}{m}.
\end{align*}

We did not expect such a simple answer  $\frac{k}{m+1}\binom{(k+1)m}{m}$ to this question about Hoppy's early adventures!

This analysis of Hoppy's early adventures covers sequences
A007226, A007228, A124724 of \cite{OEIS}, which references to \cite{AHS}.

 \bibliographystyle{plain}

\end{document}